\documentclass[12pt]{article}
\usepackage{graphicx} 
\usepackage{amsmath}
\usepackage{amssymb}
\usepackage{amsthm}
\usepackage{latexsym}

\newcommand{\comment}[1]{}
\newtheorem{theorem}{Theorem}        
\newtheorem{lemma}{Lemma}[section]
\newtheorem{remark}{Remark}[section]

\newtheorem{proposition}{Proposition}[section]

\newtheorem{definition}{Definition}[section]

\newtheorem{corollary}{Corollary}[section]
\begin{document}

\date{}

\title{{\LARGE\sf Clusters and Recurrence in 
the Two-Dimensional Zero-Temperature Stochastic
Ising Model}}

\author{
{\bf Federico Camia}\\
%\thanks{Research partially supported by the
%U.S. NSF under grant
%DMS-98-02310.}\\
{\small \tt federico.camia\,@\,physics.nyu.edu}\\
{\small \sl Department 
of Physics, New York University, New York, NY 10003, USA}\\
\and
{\bf Emilio De Santis}\\
%\thanks{Research partially supported by the
%Italian CNR.}\\
{\small \tt desantis\,@\,mat.uniroma1.it}\\
{\small \sl Universit\`{a} di Roma {\it La Sapienza\/},
Dipartimento di Matematica ``Guido Castelnuovo''}\\
{\small \sl Piazzale Aldo Moro 2, 00185 Roma, Italia}
\and
{\bf Charles M.~Newman}\\
%\thanks{Research partially supported by the
%U.S. NSF under grant
%DMS-98-03267.}\\
{\small \tt newman\,@\,courant.nyu.edu}\\
{\small \sl Courant Inst.~of Mathematical Sciences, 
New York University, New York, NY 10012, USA}\\
}

\maketitle
\begin{abstract} 
We analyze clustering and (local) recurrence of a standard
Markov process model of spatial domain coarsening. 
The continuous time process, whose state space consists of
assignments of $+1$ or $-1$ to each site in ${\bf Z}^2$, is the
zero-temperature limit of the stochastic homogeneous Ising ferromagnet
(with Glauber dynamics): the initial state is chosen uniformly at random
and then each site, at rate one, polls its $4$ neighbors and makes sure it
agrees with the majority, or tosses a fair coin in case of a tie.  
Among the main results (almost sure, with respect to both the process
and initial state) are: clusters 
(maximal domains of constant sign) are finite for times $t< \infty$, but
the cluster of a fixed site diverges (in diameter)
as $t \to \infty$; each of the two constant states is (positive) recurrent. 
%The set of nonconstant absorbing states is recurrent and 
%the set of non-absorbing states is null recurrent.
%Conjectures and scenarios for the
%positive and null recurrent states are discussed.
We also present other results and conjectures concerning positive and
null recurrence and the role of absorbing states.

\end{abstract}



Abbreviated title: Recurrence in the Stochastic Ising Model.

%\newpage

\section{Synopsis} \label{synopsis}

Consider the following Markov process, whose state
$\sigma^t$ at (continuous) time $t$ is an assignment to each site 
in ${\bf Z}^2$ of
$+1$ or $-1$. The initial state is chosen uniformly 
at random and then with rate one each site changes its value 
(resp., determines its value
by a fair coin toss) 
if it disagrees with 
three or four (resp., exactly two) of its four nearest neighbors.
This process has been much studied in the physics literature as a model of 
``domain coarsening'' (see, e.g., \cite{Bray}): clusters of constant sign 
(either $+1$ or $-1$) shrink or grow
or split or coalesce as their boundaries evolve. 
%The absorbing states are those in which every site
%agrees with at least three of its neighbors.
A more detailed
definition along with some physical motivation will be given in  
Section~2 below.

One focus of this paper is the 
study of the asymptotic growth of clusters. Let $R_* (t)$
(resp., $R^* (t)$) denote the Euclidean distance from the origin
to the closest (resp., farthest) site in its cluster that is next to the
cluster boundary (i.e., that has a neighbor of opposite
sign). It was previously proved \cite{NNS} that almost surely, 
each site flips (i.e., changes its value) 
infinitely often and thus $\liminf_{t \to \infty} R_* (t) =0$.
Among our main results are the following [with references in brackets to  
later in the paper]:
\begin{itemize}
\item  For any
$t$, almost surely, $R^* (t) < \infty$;
i.e., there is no percolation at time $t$ [Proposition~\ref{P2}].
\item $R^* (t) \to \infty$ almost surely [Proposition~\ref{P1}].
\item $\limsup_{t \to \infty} R_* (t) = \infty$ almost surely [Corollary~\ref{C1}].
%\item The rate of flips at the origin tends to zero, in probability; 
%i.e., the state is locally absorbing, with probability
%approaching one [Theorem~\ref{T2}].
\end{itemize}

The last of the three results just mentioned is a corollary of our other main
focus --- the analysis of (local) recurrence of states
$\sigma$ or (measurable) subsets ${\cal M}$ of states. We say that
$\sigma$ is \emph{recurrent} (for a given $\omega$ in the underlying
probability space) if $\sigma ^{t_k} \to \sigma$ along some
subsequence $t_k \to \infty$. We say that ${\cal M}$ is recurrent if
some $\sigma \in {\cal M}$ is recurrent.
(Related notions of recurrence for interacting particle systems
are studied in a recent paper by Cox and Klenke \cite{CK}.)
A non-recurrent $\sigma$ or ${\cal M}$
will naturally be called \emph{transient}.
%almost surely
%$\sigma^t$ agrees locally with $\sigma$ (resp., with some 
%$\sigma \in {\cal M}$) for an unbounded set of times.
Although it has not yet been proved that the 
probability distribution $\mu^t$ of
$\sigma^t$ has a unique limit as $t \to \infty$, 
nevertheless, we will classify 
a recurrent $\sigma$ (or ${\cal M}$) 
as \emph{positive recurrent} if some subsequence limit $\mu$ of $\mu^t$
(these always exist by compactness)
has $\mu (\{\sigma\}) > 0$ (or $\mu ({\cal M}) > 0$);
%the overall probability for
%$\sigma^t$ locally agreeing with $\sigma$ (resp., with some 
%$\sigma \in {\cal M}$) does not tend to zero as $t \to \infty$;
otherwise it is classified 
as \emph{null recurrent}. 

The formulation of our recurrence results involves
the absorbing states of the process, i.e., those states in which every site
agrees with at least three of its neighbors. Note that these
need not be recurrent, since our definition of recurrence is with respect
to a uniformly random initial state. It is easy to see that besides the two
constant states (identically $+1$ or identically $-1$), the absorbing states
are those whose clusters are all either half spaces or infinite strips,
and their cluster boundaries are all doubly-infinite flat lines
(either all vertical or all horizontal) separated from each other by distance
at least two. We prove the following, where $\mu$ denotes \emph{any} subsequence
limit of $\mu^t$:
\begin{itemize}
\item $\mu (\{ \hbox{non-absorbing states} \})\,=\,0$
[Theorem~\ref{T2}] or
equivalently, the rate of flips at the origin tends to zero, in probability.
Thus, almost surely, the set of non-absorbing states is not
positive recurrent.
\item Almost surely, each of the two constant states 
is positive recurrent [Theorems~\ref{T3},~\ref{T4}]. 
\end{itemize}
A natural conjecture (see Remark~\ref{R1})
is that $\mu (\{ \hbox{non-constant states} \})\,=\,0$
(or equivalently, that the 
density of cluster
boundaries tends to zero, or equivalently that $R_*(t) \to \infty$ in
probability) and thus
that, almost surely, the set of non-constant states is not
positive recurrent. 

We also prove:
\begin{itemize}
\item Almost surely, 
the set of non-constant absorbing states is recurrent [Theorem~\ref{T4}].
\item Almost surely, the 
set of non-absorbing states is null recurrent [Remark~\ref{R2}].
\end{itemize}

\noindent The non-absorbing states of the last-mentioned result 
may be restricted to those
with a cluster boundary that is flat except for a
single step of unit size next to the origin (providing the origin
with a tie among its four neighbors). There are two plausible scenarios
concerning the exact class of null recurrent states. To explain these,
we first note that cluster boundaries in recurrent states must be doubly infinite
and {\it monotonic\/}, i.e., with every finite segment either flat or
else having a Southwest--Northeast (resp., Northwest--Southeast) orientation
[Lemma~\ref{L2}]. If there were more than a single such ``domain wall''
in a recurrent state, there would be further restrictions concerning their relative
locations (see, e.g., the proof of Lemma \ref{L5}). But we conjecture that this
does not occur; indeed, we expect that one of the following
two possibilities occurs.
\begin{itemize}
\item Scenario 1: Almost surely, the null recurrent 
states are exactly all states with
single infinite monotonic domain walls while the set of all other nonconstant
states is transient.
\item Scenario 2: Almost surely, 
the null recurrent states are exactly all states with
single infinite domain walls that are either completely flat or else
have a single step, with unit size, while the set of all other nonconstant
states is transient. 
\end{itemize}

\section{Introduction} \label{introduction}

The behavior of different kinds of magnetic systems following a deep quench 
is a central topic in the study of their nonequilibrium dynamics. 
Physically, a deep quench is when a system that has reached equilibrium at
some high temperature $T_1$ has its temperature rapidly reduced to a much
lower $T_2$.
In this paper, as in much of the theoretical physics literature, we take 
$T_1 = \infty$ and $T_2 = 0$.
Rigorous
and nonrigorous results have been obtained on different questions that arise 
naturally in this context, such as the formation of domains, their subsequent 
evolution, spatial and temporal scaling properties and related problems (for a
review, see \cite{Bray}).
In particular, in the context of zero-temperature, stochastic Ising models
with nearest-neighbor interactions, the question of whether the spin configuration
eventually settles down to a final state has been 
addressed rigorously and answered
for a number of different models \cite{GNS, NNS}. 
A closely related issue is that of persistence, concerning the fraction of
sites that have not flipped at all by time $t$, and its asymptotic behavior
as $t \to \infty$ \cite{Derrida,DHP,NS2,Stauffer}.

In this paper, we consider the 
zero-temperature, stochastic (homogeneous) Ising model $\sigma^t$
on ${\bf Z}^d$ with nearest-neighbor ferromagnetic interactions.
In dimension one, the model is the same as the $d=1$ voter model,
which is well understood (see, e.g., Chapter 2 of \cite{Durrett} or
Chapter V of \cite{Liggett} and references therein).
We study in detail the case $d=2$, for which it is known
that there is not a unique limiting state \cite{NNS}, 
and consider questions about
the limits along subsequences of time and 
the nature of clusters of sites of the same sign. 
The aim of the paper is to give a 
picture of the system for very long times, showing
what kinds of events and 
states are seen as $t \to \infty$ and
what are instead ``forbidden''.
These are basically questions of recurrence and transience.

In the remainder of this section, we define the model precisely
and discuss some results and open problems about this and related models. 
The stochastic process $\sigma^t = \sigma^t(\omega)$ corresponds to the
zero-temperature limit of Glauber dynamics for an Ising model with 
(formal) Hamiltonian
\begin{equation} \label{Hamiltonian}
{\cal H} = - \sum_{ \stackrel{ \{x,y \} }{||x-y||=1} } 
J_{x,y} \sigma_x \sigma_y .
\end{equation}
Here the state space is
$ {\cal S} = \{ -1, +1 \}^{{\bf Z}^d} $, the space of
(infinite-volume) spin configurations $\sigma$, 
and $\Vert \cdot \Vert$ denotes Euclidean 
length. The initial spin 
configuration $\sigma^0$ is chosen from a symmetric Bernoulli
product measure (denoted $\mu^0$), 
corresponding physically to a deep quench from infinite
temperature.
(We note that the case of an {\it asymmetric\/} initial
$\mu^0$ has also been studied, both on ${\bf Z}^d$ \cite{FFS}
and on other lattices \cite{Howard}.)
The continuous time dynamics is defined by means 
of independent (rate $1$) Poisson processes
at each site $x$ corresponding to those 
times (which we think of as ``clock rings'') when a spin flip 
$(\sigma_x^{t+0} = - \sigma_x^{t-0})$ is \emph{considered}. 
If the change in energy,
\begin{equation}
\Delta {\cal H}_x (\sigma) = 
2 \sum_{y: \Vert x-y \Vert = 1} J_{x,y} \sigma_x \sigma_y ,
\end{equation}
is negative (or zero or positive), 
then the flip is done with probability $1$ (or $1/2$
or $0$).
We think of associating a fair coin toss to each
clock ring, which we use as a tie-breaker only when
$\Delta {\cal H}_x (\sigma) = 0$.
Let us denote by $P_{dyn}$ the probability measure for the
dynamics realizations of clock rings and tie-breaking coin
tosses and then denote by $P = \mu^0 \times P_{dyn}$
the joint probability 
measure on the space $\Omega$ of the initial
configurations $\sigma^0$ and  
realizations of the dynamics. An element of 
$\Omega$ will be denoted by $\omega$. 

The model that we will study in 
this paper is the homogeneous ferromagnet, where
$J_{x,y} \equiv 1$ for all $ \{ x,y \} $. 
In this model, when the clock at 
site $x$ rings, $\sigma_x$ flips with 
probability $1/2$ if it disagrees with exactly 
$d$ neighbors and with probability $1$ if 
it disagrees with more than $d$ neighbors; it does not
flip if it disagrees with less than $d$ neighbors.
In the first case, the spin flip leaves the energy unchanged, 
$\Delta {\cal H}_x (\sigma) = 0$, while in 
the second case the spin flip lowers the
energy, $\Delta {\cal H}_x (\sigma) < 0$.
A very useful result of Nanda, Newman 
and Stein (see Theorem 3 and the following
remark in \cite{NNS})
 states that the number of energy 
lowering spin flips at any site is almost surely finite. 

Disordered models, in which a realization 
${\cal J}$ of the $J_{x,y}$'s is chosen from the
(independent) product measure  
of some probability measure $\nu$, are also
studied in the literature (see, e.g., \cite{NNS, NS}) 
and have different properties, but we will 
not deal with those models here. For
the homogeneous ferromagnet in dimensions $d=1$ \cite{Arratia,CG}
and $d=2$ \cite{NNS}, 
$\sigma^{\infty}(\omega) = \lim_{t \to \infty} \sigma^t (\omega)$ 
does not exist; indeed, for 
almost every $\omega$ and for every 
$x$, $\sigma_x^t(\omega)$ flips infinitely
many times. 
For $d>2$, little is known rigorously, but 
numerical studies \cite{Stauffer} suggest
that the same is true up to $d = 4$, while
$\sigma^{\infty}(\omega)$ might perhaps exist for $d > 4$. 

When the limit does not exist, it is natural to ask what happens
to the measure describing the state of the 
system as $t \to \infty$. A natural way 
to approach this question is by looking 
at the clusters of sites with the same sign and at the
domain walls between such clusters. 
The two descriptions are basically equivalent and we will 
use both, depending on the type of problem. 

It is a direct consequence of a result of 
Harris \cite{Harris, Liggett} that the distribution $\mu^t$ of
$\sigma^t$ satisfies the FKG property for any time $t$. 
In dimension $d=2$, we will use this and a 
result of Gandolfi, Keane and Russo \cite{GKR} to show that
at any time $t$, neither $+1$ nor $-1$ spins percolate;
i.e., the clusters are almost surely finite. 
We will also show, however, that the diameter
of the cluster at the origin  
almost surely diverges as $t \to \infty$, for any $d$. 
In dimension two, it is a natural conjecture that for 
large enough times the system will be, with 
$P$-probability close to $1$, locally in a 
$+1$ or in a $-1$ phase or, equivalently, that the 
density of domain walls tends to zero.
In fact, it is expected that the density is of order $t^{-1/2}$ as
$t \to \infty$ (see, e.g., \cite{Bray}).
Although we are not able to prove this, Theorem \ref{T3} in Section~5
below points in that direction.

If the above conjecture is true, then, by symmetry, it automatically follows
that the distribution $\mu^t$ of $\sigma^t$ has the unique limit, as
$t \to \infty$, of $\frac{1}{2} \delta_{+1} + \frac{1}{2} \delta_{-1}$, where
$\delta_{\eta}$ is the probability measure assigning probability one to the
constant ($\equiv \eta$) spin configuration.
But $\mu^t$ is the overall distribution of $\sigma^t$, taking into
account that the initial state is random and distributed by the Bernoulli
product measure $\mu^0$.
If instead, we condition on $\sigma^0$ and consider the conditional distribution
$\mu^t[\sigma^0]$ (for almost every $\sigma^0$), it is unclear whether
that should still converge as $t \to \infty$ to
$\frac{1}{2} \delta_{+1} + \frac{1}{2} \delta_{-1}$, or rather there should
be multiple subsequence limits (presumably all of the form
$\alpha \delta_{+1} + (1-\alpha) \delta_{-1}$) along different 
$\sigma^0$-dependent subsequences of time.
The latter situation would be an example of ``Chaotic Time Dependence'' (CTD)
\cite{FIN} (see also \cite{NS1}).
CTD is known not to occur for the $d=1$ version of our model (equivalent to the
voter model), but has been proved to occur in a disordered $d=1$ voter
model \cite{FIN, FIN1}.

\section{Percolation results} 

In this section we present two propositions about clusters of constant sign.
For every $x \in {\bf Z}^d$, let us denote by $C_x(t) = C_x(\sigma^t)$ the 
\emph{cluster at site $x$ at time $t$}. 
$C_x(t)$ is defined as the maximal subset of ${\bf Z}^d$
satisfying the following properties:
\begin{itemize}
\item $x \in C_x(t)$, 
\item $C_x(t)$ is connected
(in the sense that if $y$ and $z$ are both in $C_x(t)$, there exists a sequence 
$\zeta_i$, $i = 0, 1, 2, \dots, n$, of sites of $C_x(t)$ with 
$||\zeta_{i+1} - \zeta_i||=1$ and with
$\zeta_0 = y$ and $\zeta_n = z$),
\item if $y,z \in C_x(t)$, then $\sigma^t_y = \sigma^t_z$.
\end{itemize}
$|C_x(t)|$ will denote the number of sites in $C_x(t)$. 
The origin is denoted by $o$ and so $C_o(t)$ is the cluster at the origin.
\\

Our first result, and the only one valid for all $d$,
concerns the growth of $C_o(t)$ with time.
It remains valid in a very general setting (see Theorem 3 and the
following remark in \cite{NNS}) and in particular applies to
our Markov process when the initial state is chosen according to any
translation-invariant measure.

\begin{proposition} \label{P1} For
any $d$, the size of the cluster at the origin diverges almost surely 
as $t \to \infty$: $\lim_{t \to \infty} |C_o(t)| = \infty$.
\end{proposition}

{\bf Proof.} We will prove the proposition by contradiction. Suppose that
the conclusion is not true; then with 
positive probability, $ \liminf_{t \to \infty} |C_o (t)| < \infty $
and so there exist $M < \infty$ and 
a sequence of times $\{ t_k \}_{k \in {\bf N}}$ 
with $t_k \to \infty$ such that $ |C_o (t_k)| < M $ for $k = 1, 2, \dots $.
Without loss of generality, we may assume that $ t_{k+1} > t_k + 1 $.
There are only finitely 
many shapes (lattice animals) that the cluster at the origin
can have at times $t_k$ when $ |C_o (t_k)| < M $. For each such lattice 
animal, there is some ordered finite sequence of clock rings and outcomes
of tie-breaking coin tosses inside a fixed finite ball that would
cause the cluster to shrink 
to a single site at the origin which would then have 
an energy lowering spin flip. It follows that for some $\delta>0$ and any 
$\sigma \in {\cal S}$ such that $|C_o(\sigma)|<M$,
\begin{equation}
P( \hbox{ origin flips at time } t \in (t_k, t_{k}+1) \hbox{ with }
\Delta {\cal H}_o < 0 \mid  \, \sigma^{t_k} = \sigma ) \geq \delta.
\end{equation}
By the Markov property of the process and our supposition that the conclusion
of the proposition is false, this would imply
that the spin at the origin $\sigma_o$ flips infinitely many times with
$ \Delta {\cal H}_o < 0 $ with positive probability, which 
contradicts a result
of Nanda, Newman and Stein \cite{NNS} mentioned in Section~\ref{introduction}
of the paper. \fbox{}
\\

We set $d=2$ now and for the rest of the paper. 
%Since, as mentioned before, in dimension two,
%the state $\sigma^t$ does not have a unique limit as $t \to \infty$
%\cite{NNS}, we are interested in the limits along subsequences of time.
%Let us introduce the set of all such possible limiting configurations,
%\begin{equation}
%{\cal W} = {\cal W} (\omega) = \{ \tilde \sigma \in {\cal S} : 
%\sigma_x^{t_k(\omega)} (\omega) \to \tilde \sigma_x \hbox{ for all }
%x \in {\bf Z}^d \hbox{ and some } \{ t_k(\omega) \}_{k \in {\bf N}}
%\hbox{ with } t_k(\omega) \uparrow \infty \}. 
%\end{equation}
Our second result is a direct consequence of a result of Harris 
\cite{Harris, Liggett} and one of Gandolfi, Keane and Russo \cite{GKR}.

\begin{proposition} \label{P2}
At any (deterministic) time,
there is no percolation of clusters of spins of the same sign:
$ P(|C_o(t)| = \infty) = 0 $ for all $t \geq 0$.
\end{proposition}

{\bf Proof.} First note that the measure $\mu^t$ describing the state
$\sigma^t$ of the system at time $t$ is invariant and
ergodic under ${\bf Z}^2$-translations. This is so  because 
the same is true for both
$\mu^0$ and $P_{dyn}$ and hence also for $P$.
Applying a result of Harris \cite{Harris, Liggett}, we also have that
$\mu^t$ satisfies the FKG property, i.e., 
increasing functions of the spin
variables are positively correlated (this follows from the FKG property of
$\mu^0$ and the attractivity of the Markov process).
Then it follows from a result of Gandolfi, Keane and Russo \cite{GKR} that if 
percolation of, say, $+1$ sites were to occur, all the $-1$ clusters 
would have to be finite. 
Because of the symmetry of the model under a global spin flip, however,
percolation of $+1$ sites with positive probability implies the same for 
$-1$ sites. 
Then, using the ergodicity of the measure, we would see simultaneous
percolation of both signs, thus obtaining a contradiction. \fbox{}
\\

\section{Preliminary recurrence results}

We now introduce the contour representation in the dual lattice 
${{\bf Z}^2}^* \equiv {\bf Z}^2 + (1/2, 1/2)$, following the notation of 
\cite{GNS}.
A (dual) site in ${{\bf Z}^2}^*$ may be identified with the plaquette $p$ 
in ${\bf Z}^2$ of which it is the center. 
The edge $ \{ x, y \}^* $ of ${{\bf Z}^2}^*$, dual to (i.e., perpendicular 
bisector of) the edge $ \{ x, y \} $ of ${\bf Z}^2$, is said to be 
\emph{unsatisfied} (with respect to a given spin configuration
$\sigma \in {\cal S}$) if $\sigma_x \not= \sigma_y$ (and \emph{satisfied}
otherwise).
Denote by $\Gamma$ the set of unsatisfied (dual) edges.
Given a finite rectangle $\Lambda$ of ${\bf Z}^2$, $\Lambda^* \subset {{\bf Z}^2}^*$ 
consists of the dual sites corresponding to the plaquettes contained
in $\Lambda$. 
$\Gamma (\Lambda^*)$ is the set of unsatisfied (dual) edges bisecting the edges
connecting sites in $\Lambda$.
Note that the outermost edges of $\Gamma (\Lambda^*)$ have one endpoint
just outside of $\Lambda^*$.
A (site self-avoiding) path in ${{\bf Z}^2}^*$ using only unsatisfied edges
will be called a \emph{domain wall}; it is simply a path along the cluster
boundaries of $\sigma$.
If $\Gamma(\Lambda^*)$ is not empty, it can contain one or more domain walls.
Since domain walls are the boundaries between clusters of sites with
different sign, they can always be extended to form a closed loop 
or a doubly infinite path.
Every $\Gamma ({{\bf Z}^2}^*)$ configuration corresponds to two spin configurations
related by a global spin flip.
The Markov process $\sigma^t$ determines a process $\Gamma^t$, that is easily 
seen to also be Markovian. The transition associated with a spin flip at
$x \in {\bf Z}^2$ is a local ``deformation'' of the contour $\Gamma^t$ at the
(dual) plaquette that contains $x$; this deformation
interchanges the satisfied and unsatisfied edges of that plaquette.
The only transitions with nonzero rates are those where the number of
unsatisfied edges starts at $k = 4$ or $3$ or $2$ and ends at $0$ or $1$ or $2$,
respectively; transitions with $k = 4$ or $3$ (resp. $2$) correspond to
energy-lowering (resp., zero-energy) flips and have rate $1$ (resp. $1/2$).
We will continue to use the terms flip, energy-lowering, etc. for the
transitions of $\Gamma^t$.

We continue with some definitions and lemmas.

\begin{definition} \label{recurrent}
Let $Z_t$ be a continuous-time Markov process with state space
${\cal Z}$ and time homogeneous transition probabilities. 
For $A$ a (measurable)
subset of ${\cal Z}$ we say that 
$A$ \emph{recurs} if $ \{ \tau : Z_{\tau} \in A \} $ 
is unbounded, and we say that $A$ is \emph{eventually absent (e-absent)}
if it recurs with zero probability. For our stochastic Ising model,
the restriction  $\sigma|_\Lambda$ of some $\sigma \in {\cal S}$
to $\Lambda \subset {\bf Z}^2$ will be called e-absent if
$\{ \sigma' \in {\cal S} :\, \sigma'|_\Lambda = \sigma|_\Lambda \}$
is e-absent. 
\end{definition}

Note that if a contour event $A$, specified by $\Gamma(\Lambda^*)$,
is e-absent, then any $\sigma|_\Lambda$ consistent with $A$ is
also e-absent. Note further that our definition in Section 1
for recurrence of $\sigma \in {\cal S}$ is that for
every finite $\Lambda$, the restriction $\sigma|_\Lambda$ recurs.
Thus almost sure transience is implied by, but not equivalent to
e-absence of $\sigma|_\Lambda$ for all large $\Lambda$ since,
a priori, it could be that $\sigma|_\Lambda$ recurs with nonzero
probability, tending to zero as $\Lambda \to {\bf Z}^2$. 

By $Q_L$ we denote the square of size $2L+1$ centered at the origin, that is
the set of all $x=(x_1, x_2) \in {\bf Z}^2$ such that $x_i \in 
\{-L, \dots, L\}$.
($Q_L(x)$ will be used later to
denote the square of size $2L+1$ centered at $x$; i.e., 
$Q_L(x) = Q_L + x$). 

\begin{lemma} \label{L4}
If $A$ is e-absent, then $P (Z_t \in A) \to 0$ as $t \to \infty$.
In particular, if we denote by $S_L$ the set of $\sigma \in {\cal S}$
such that $\sigma|_{Q_L}$ is e-absent, then
\begin{equation}
\lim_{t \to \infty} P(\sigma^t \in S_L) = 0 .
\end{equation}
\end{lemma}

{\bf Proof.}
Suppose the lemma were false. Then there would exist $\delta > 0$ and
a sequence of times $t_k \uparrow \infty$ such that for all $k$,
$ P(Z_{t_k} \in A) > \delta $. But then it would follow that
\begin{equation} \label{impossible}
P(Z_t \hbox{ recurs } ) > \delta ,
\end{equation}
contradicting the e-absence of $A$, as a consequence of the standard
fact that
\begin{equation} \label{lemma}
 P(B_k) > \varepsilon \hbox{ for all } k \in {\bf N} \hbox{ implies }
P( B_k \hbox{ occurs infinitely often } ) > \varepsilon. ~\fbox{}
\end{equation}

The following lemma is essentially the same as Lemma 8 of \cite{GNS}
(where a more detailed proof may be found). 
We say that a domain wall in ${{\bf Z}^2}^*$ 
is \emph{monotonic} if, for one of the two directed path versions of the
domain wall, either every move is to the North or East or else
every move is to the South or East. 

\begin{lemma} \label{L2}
The event $\{ \Gamma^t(Q_L^*) \hbox{ contains a non-monotonic domain wall} \}$ 
is e-absent.
%\begin{equation}
%P( \Gamma^t(Q_L^*) \hbox{ contains a non-monotonic domain wall} \,\, i.o. ) 
%= 0.
%\end{equation}
\end{lemma}

{\bf Proof.} 
A non-monotonic domain wall in $\Gamma(Q_L^*)$ can always be modified through 
local deformations (corresponding to appropriate spin flips of sites in $Q_L$)
to give a contour configuration $\Gamma'(Q_L^*)$ with three 
(or four) domain wall
edges surrounding some plaquette of $Q_L^*$. The corresponding
spin configurations then have a site that disagrees 
with three (or four) neighbors, which can undergo 
an energy-lowering spin flip. As in the proof of Proposition \ref{P1},
the existence of such local deformations (or sequence of spin flips)
means that there is a bounded away from
zero probability (corresponding to an appropriate 
sequence of clock rings and tie-breaking coin tosses) of an 
energy-lowering spin flip during the next unit time interval.
If the claim were not true, then the event that $\Gamma(Q_L^*)$ contains a 
non-monotonic domain wall would recur and there would be 
a nonzero probability of
infinitely many energy lowering spin flips in $Q_L$, which 
would contradict an 
already mentioned theorem of Nanda, Newman and Stein \cite{NNS}. \fbox{}
\\

The next lemma provides a geometric upper bound which is one of the key
technical results of this paper.
In particular, it will be used in the proof of Theorem \ref{T1}
below.
Before we can state the lemma, we need a definition.
For 
%$x \in {\bf Z}^2$, $l \in {\bf N}$ and 
%$t \geq 0$, 
a given $\sigma \in {\cal S}$,
let $M_L(\sigma)$ denote the total number
of corners in $\Gamma(Q_L^*)$, i.e., pairs of perpendicular edges that
meet at a site in $Q_L^*$.
$S_L^{\,c}$, the complement
of $S_L$, is the set of $\sigma$'s such that $\sigma|_{Q_L}$
is not e-absent.

\begin{lemma}  \label{L5} For $\sigma \in S_L^{\,c}$, the number
$M_L(\sigma)$ of corners is bounded by
%For any spin configuration $\sigma$ such that $\sigma|_{Q_L}$ is not e-absent
%($\sigma \in S_L^{\,c}$), the number of sites $x \in Q_L$ such that 
%$I_{l,x}(\sigma) = 1$ has the following upper bound: 
%
\begin{equation} \label{corners}
\max_{\sigma \in S_L^{\,c}} M_L(\sigma) \leq 4 (2L+1).
%\max_{\sigma \in S_L^{\,c}} \sum_{x \in Q_L} I_{l,x}(\sigma)  \leq  
%[ (16\,l^2 + 4\,l) (2L+1) 
%\wedge (2L+1)^2 ] .
\end{equation}
\end{lemma}

{\bf Proof.} 
%For a given $\sigma$, let $M_L$ denote the total number
%of corners in $\Gamma(Q_L^*)$, i.e., pairs of perpendicular edges that 
%meet at a site in $Q_L^*$. It is an elementary fact that if 
%$I_{l,x}=1$, then there must be a 
%``nearby'' corner, i.e., one in $\Gamma(Q_l^*(x))$.
%This corner will be in $\Gamma(Q_L^*)$ as long as $x$ is in $Q_L$ and
%its distance from $Q_L^c$ exceeds $l$. Thus each corner in $\Gamma(Q_L^*)$
%allows (at most) $4l^2$ sites $x$ to have $I_{l,x}=1$, and combining
%these with sites too close to $Q_L^c$ yields the bound,
%\begin{equation} \label{corners2}
%\sum_{x \in Q_L} I_{l,x}(\sigma) \leq
%M_L \cdot 4L^2+\,l\cdot 4 (2L+1) .
%\end{equation}
%The remainder of the proof is a demonstration that for $\sigma \in S_L^{\,c}$,
%$M_L \leq 4(2L+1)$, which immediately yields (\ref{corners}).
%
By Lemma \ref{L2}, any 
$\Gamma(Q_L^*)$ that is not e-absent 
may be partitioned into edge-disjoint monotonic domain walls,
$\gamma_1, \dots, \gamma_m$, whose endpoints are just outside of
$Q_L^*$. Let us denote by $M(\gamma_i)$ the number of corners
in $\gamma_i$. Clearly, if $m=1$, then $M_L = M(\gamma_1) \leq 2(2L+1)$.
When $m > 1$, we need to consider the geometric constraints on the
$\gamma_i$'s required for $\sigma \in S_L^{\,c}$. First, there are the special
$m=2$ cases where $\Gamma(Q_L^*)$ is a ``cross,'' i.e., the union
of one flat horizontal and one flat vertical domain wall; here $M_L = 4
\leq 4(2L+1)$.
We claim that in all remaining cases, the $\gamma_i$'s
are site-disjoint and in fact, 
their spanning rectangles $R(\gamma_i)$, defined so
that two of the vertices of $R(\gamma_i)$
are the endpoints of 
$\gamma_i$, must also be site-disjoint. 

Note that, depending on
whether $\gamma_i$ connects two opposite or two adjacent sides of the
boundary of $Q_L^*$, the rectangle $R(\gamma_i)$ can be classified
as either vertical or horizontal or as one of four corner-types.
The reason the $R(\gamma_i)$'s must be site-disjoint (except in
a cross configuration) is that otherwise there would
exist a sequence
of spin flips (corresponding to a sequence of clock-rings and  
tie-breaking coin tosses) in $Q_L$ that would deform $\Gamma(Q_L^*)$
into a contour configuration with a non-monotonic domain wall
(see Figure 1 and Lemmas 9 and 10(i) of \cite{GNS}). 
As in the proof of Lemma \ref{L2}, this shows that intersecting
$R(\gamma_i)$'s must be e-absent.

\begin{figure}[!ht]
\begin{center}
\includegraphics[width=10cm]{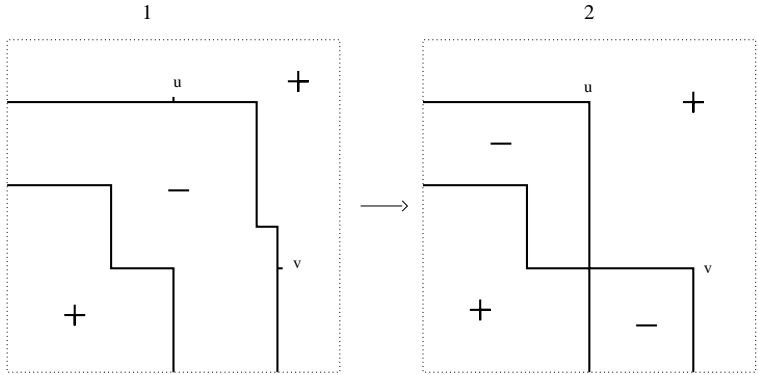}
\caption{An example in which
configuration 1 is deformed to configuration 2
by spin flips with $\Delta {\cal H} = 0$.
Configuration 2 has non-monotonic domain walls.}\label{f1}
\label{int}
\end{center}
\end{figure}

When the  $\gamma_i$'s are site-disjoint, $M_L = \sum_i M(\gamma_i)$,
and if
$R(\gamma_i)$ has sides of length $l_1^i$ and $l_2^i$, then $M(\gamma_i) \leq
2 \min (l_1^i, l_2^i)$. 
Now it is easily seen that the sum of the shorter sides
of these non-overlapping $R(\gamma_i)$'s is bounded by
twice the linear dimension of $Q_L$ (i.e.
$2(2L+1)$), yielding $M_L \leq 4(2L+1)$ as claimed; indeed the
worst case is when $m=4$ and each $\gamma_i$ is of corner-type,
with the $R(\gamma_i)$'s almost overlapping. \fbox{}
\\

One of the main results of the paper,
Theorem \ref{T2} below, concerns the probability
that, for large $t$, $\sigma^t$ is locally in an absorbing state.
As noted in Section 1, the absorbing states are constant either on infinite
horizontal lines or infinite vertical lines.
The next result of this section is a more technical theorem about the
density of domain wall corners, from which
Theorem \ref{T2} follows easily. For $t\geq 0$, 
and $x^* \in {{\bf Z}^2}^*$, we define 
$F_{x^*}(t)$ to be the event that there is a corner in
$\Gamma ^t ({{\bf Z}^2}^*)$ with vertex $x^*$.
%$\sigma^t$ is constant on the horizontal
%and/or on the vertical line segment of length $2l+1$ centered at $x$.
%We also write $F_l(t)$ for $F_{l,o}(t)$.

\begin{theorem} \label{T1}
For any $x^* \in {{\bf Z}^2}^*$, $ \lim_{t \to \infty } P (F_{x^*} (t)) = 0 $. 
\end{theorem}

{\bf Proof.} 
%Define $N_L(t)$ to be the number of sites $x$ in
%the square $Q_L$ such that $F_{l,x} (t)$ does not occur. This
%is the same as the left-hand side of (\ref{corners2}) with 
%$I_{l,x} (\sigma)$ replaced by $1_{(F_{l,x} (t))^c} = I_{l,x} (\sigma^t)$.
If $F_{x^*}(t)$ occurs, then there is at least one (and at most four)
corners at $x^*$.
Thus $\tilde{M}_L(\sigma^t)$, the number of $x^*$'s in $Q_L^*$ such that
$F_{x^*}(t)$ occurs, is bounded by the number $M_L(\sigma^t)$ of corners,
and so by translation invariance and elementary arguments,

\begin{eqnarray}
P( F_{x^*} (t) ) & = & E(1_{F_{x^*} (t)}) = \frac{1}{(2L)^2} E(\tilde{M}_L(\sigma^t))
\nonumber \\
{} & = & \frac{1}{(2L)^2} \left[ P(\sigma^t \in S_L) E(\tilde{M}_L(\sigma^t) |
\sigma^t \in S_L) + P(\sigma^t \in S_L^c) E(\tilde{M}_L(\sigma^t) | \sigma^t \in S_L^c)
\right]
\nonumber \\
{} & \leq & \frac{1}{(2L)^2}  \left [ P(\sigma^t \in S_L) \cdot
(2L)^2 + \, E(M_L(\sigma^t)\,|\, \sigma^t \in S_L^c)  \right ] \nonumber \\
{} & \leq & P(\sigma^t \in S_L) \, + \, \frac{1}{(2L)^2}
\max_{\sigma \in S_L^c} M_L(\sigma) \, .
%\sum_{x \in Q_L} I_{l,x} (\sigma) \, .
\end{eqnarray}
The proof is completed by using Lemmas \ref{L4}
and \ref{L5} to observe that the two 
terms in the final expression 
can be made small by appropriate choice of large $L$ and $t$. \fbox{}
\\

\section{Main recurrence results}

Let us now define two sets of states 
that play an important role in recurrence:
\begin{equation}
{\cal C}  = \{ \sigma \in {\cal S} \, : 
\, \sigma \equiv 1  \hbox{ or }  \sigma 
\equiv -1 \},
\end{equation}
\noindent
the set of constant spin configurations, and
\begin{equation}
{\cal A}  =   
\{ \hbox{ absorbing states }   \}, 
\end{equation}
\noindent
consisting of the two constant spin configurations 
in ${\cal C}$ together with all spin configurations
corresponding to contour configurations that are unions of
doubly infinite flat domain walls (either all horizontal or all
vertical)
separated by at least two lattice spacings. These are the only absorbing 
states in $d=2$, since only for these does 
every spin agree with a strict majority of its neighbors.
%\begin{equation}
%{\cal G}  ={\cal C} \cup  
%\{ \hbox{ single monotonic domain wall states }   \}, 
%\end{equation}
%\begin{equation}
%{\cal G^*}  ={\cal C} \cup  
%\{ \hbox{ single flat domain wall states }   \}.
%\end{equation}

We also introduce the set of all spin configurations 
${\tilde \sigma}$ that agree with a given subset of ${\cal S}$
inside the square $Q_L$: for $ {\cal U}  \subset {\cal S}$, define
\begin{equation}
{\cal U}_L  = \{ \tilde \sigma \in {\cal S}  :
\exists \, \sigma \in  {\cal U}
\hbox{ with } {\sigma }|_{Q_L} = {\tilde \sigma }|_{Q_L} \}.
\end{equation}

%\begin{definition} Given $ {\cal U}  \subset \{-1 , 1 \}^{ {\bf Z}^d  } $ 
%and $L < \infty$, define 
%\begin{equation} 
%{\cal U}_L  = \{ \tilde \sigma \in \{-1 , 1 \}^{ {\bf Z}^d } : 
%\exists \, \sigma \in  {\cal U} 
%\hbox{ with } {\sigma }|_{Q_L} = {\tilde \sigma }|_{Q_L} \}. 
%\end{equation} 
%\end{definition}

The next theorem states that $\sigma^t$, in any finite region,
for large enough $t$, agrees
with some absorbing state with probability 
arbitrarily close to $1$. It is straightforward to see that this
is equivalent to saying that the 
rate of flips 
at the origin goes to zero in probability as $t \to \infty$.

\begin{theorem} \label{T2}
For all $L \in {\bf N}$,
$\lim_{t \to \infty} P( \sigma^t \in {\cal A}_L ) = 1.$
\end{theorem}

{\bf Proof. } We will prove the theorem by contradiction. 
If the claim is not true, there exist $L > 0$,
$\delta > 0$ and a sequence $t_k$ with $t_k \to \infty$ such that for all $k$,
$ P( \sigma^{t_k} \notin {\cal A}_L ) > \delta $.
Now, the event $ \{ \sigma^t \notin {\cal A}_L \} $ corresponds to
either having 
\begin{itemize}
\item at least one non-flat domain wall inside $Q_L$, or else
\item two flat domain walls at distance one apart.
\end{itemize}
The first case corresponds to having at least one site $x^* \in Q_L^*$ for
which the event $F_{x^*}(t)$, that there is a corner at $x^*$,
occurs.
%for all $l>0$ (every site
%corresponding to a corner made by a domain wall has this property).
The second case corresponds to a configuration that can become non-monotonic
with one flip and is therefore e-absent.
Thus
\begin{equation}
P(\sigma^t \notin {\cal A}_L) \leq 
\sum_{x^* \in Q_L^*} P( F_{x^*}(t) ) + P(\sigma^t \in S_L),
\end{equation}
and the right-hand side goes to zero as $t \to \infty$ by Theorem
\ref{T1} and Lemma \ref{L4}. This
completes the proof. \fbox{} \\

The next theorem, combined with the first part of Theorem \ref{T4},
shows that each of the two
constant states is positive recurrent in the sense of Section 1.

\begin{theorem} \label{T3}
Let $C^+_L (t)$ (resp., $C^-_L (t)$) denote the event that $\sigma^t$
is constant and equal to $+1$ (resp., to $-1$) on the square $Q_L$.
Then for any $L < \infty$, 
\begin{equation}
\liminf_{t \to \infty} P(C^+_L (t))\,
= \, \liminf_{t \to \infty} P(C^-_L (t)) \, \geq  \, 1/4. 
\end{equation}
%For all finite $L$,
%$ \liminf_{t \to \infty} P( \sigma^t \in {\cal C}_L ) \geq 1/2 $.
\end{theorem}

{\bf Proof.} 
%Given any square centered at the origin $Q_L$, we 
We introduce the
following events that, like $C^+_L (t)$, are increasing in the FKG sense: 
\begin{eqnarray}
%C^+_L (t) & = & \{ \omega : \sigma^t_x (\omega) \equiv +1, 
%\forall x \in Q_L \}, \\
V^+_L (t) & = & \{ \omega : \sigma^t|_{Q_L} (\omega) 
\hbox{ has a vertical line of } 2L +1 \hbox{ sites that are } +1 \}, \\
H^+_L (t) & = & \{ \omega : \sigma^t|_{Q_L} (\omega) 
\hbox{ has a horizontal line of } 2L+1 \hbox{ sites that are } +1 \}.
\end{eqnarray}
Note that $C^+_L (t) \subset V^+_L (t), H^+_L (t)$.
We also define the corresponding events with $+$ replaced by $-$.
With these definitions, we have
\begin{equation}
V^+_L (t) \cap H^+_L (t) \subset C^+_L (t)
\cup \{ \sigma^t \notin {\cal A}_L \} 
\end{equation}
and therefore
\begin{equation}
P(V^+_L (t) \cap H^+_L (t)) \leq P(C^+_L (t))
+ P( \sigma^t \notin {\cal A}_L ) .
\end{equation}
Using the fact that $V^+_L (t)$ and $H^+_L (t)$ are increasing
events and the FKG property of the distribution of $\sigma^t$
(see the proof of Proposition \ref{P2}), we get
\begin{equation} \label{fkg}
P(C^+_L (t)) \geq P(V^+_L (t)) P(H^+_L (t))
- P( \sigma^t \notin {\cal A}_L ) .
\end{equation}

Because of the ``striped'' nature of the absorbing states,
\begin{equation} \label{absstates}
P(\sigma^t \in {\cal A}_L) = P(\{ \sigma^t \in {\cal A}_L \} \cap H^+_L (t)) 
+ P(\{ \sigma^t \in {\cal A}_L \} \cap V^-_L (t)).
\end{equation}
By the symmetries of the model, the two terms in the right hand side of
(\ref{absstates}) are equal and therefore
\begin{equation}
P(\sigma^t \in {\cal A}_L) = 2 P(\{ \sigma^t \in {\cal A}_L \} \cap H^+_L (t)).
\end{equation}
Thus
\begin{eqnarray}
P(H^+_L(t)) & = & P(\{ \sigma^t \in {\cal A}_L \} \cap H^+_L (t))
+ P(\{ \sigma^t \in {\cal A}_L \}^c \cap H^+_L (t)) \\
& = & \frac{1}{2} P(\sigma^t \in {\cal A}_L) 
+ P(\{ \sigma^t \in {\cal A}_L \}^c \cap H^+_L (t)).
\end{eqnarray}
Applying Theorem \ref{T2} and symmetry, we obtain
\begin{equation}
\lim_{t \to \infty} P(V^+_L (t)) =
\lim_{t \to \infty} P(H^+_L (t)) = 1/2.
\end{equation}
Taking the $\liminf$ of both sides of (\ref{fkg}) and using Theorem \ref{T2}
once more, we have
\begin{equation}
\liminf_{t \to \infty} P(C^+_L (t)) \geq
\lim_{t \to \infty} P(V^+_L (t)) P(H^+_L (t)) = 1/4.~\fbox{}
\end{equation}
%Therefore,
%\begin{equation}
%\liminf_{t \to \infty} P(\sigma^t \in {\cal C}_L) =
%\liminf_{t \to \infty} P(C^+_L (t) \cup C^-_L (t))
%\geq 1/2 .~\fbox{}
%\end{equation}

\begin{remark} \label{R1} A natural conjecture is
that the system is in a constant ($+1$ or $-1$) state
with probability approaching $1$ as $t \to \infty$, i.e., for all $L$,
\begin{equation} \label{conjecture}
\lim_{t \to \infty} P( \sigma^t \in {\cal C}_L ) = 1.
\end{equation}
This is equivalent to the
conjecture for our $d=2$ process that \emph{clustering}
occurs:
\begin{equation}
P(\sigma^t_x \neq \sigma^t_y) \to 0 
\hbox{ as } t \to \infty \hbox{ for any } x,y
\in
{\bf Z}^2 .
\end{equation}
In $d=1$, our process is the 
same as the one-dimensional voter model, for which 
clustering is
known to occur \cite{HH} (see also, e.g., \cite{Durrett, Liggett}).
\end{remark}

The next result is a corollary of Theorem \ref{T3}. Recall that $R_*(t)$
denotes the Euclidean distance from the origin to the closest
site in its cluster that is next to the cluster boundary. More precisely,
given a subset $\Lambda$ of ${\bf Z}^2$ and $x \in {\bf Z}^2$, define
the distance $ d(\Lambda, x) = \inf_{y \in \Lambda} ||x-y|| $.
The inner boundary of $\Lambda$ is $ \partial \Lambda = \{ x \in \Lambda :
\exists y \notin \Lambda \hbox{ with } ||x-y|| = 1 \} $.
Then $R_* (t) = d( \partial C_o(t), o )$.

\begin{corollary} \label{C1}
$ \limsup_{t \to \infty} R_*(t) = \infty $ almost surely.
\end{corollary}

{\bf Proof.} Define the event
\begin{equation} 
A_L = \{ \limsup_{t \to \infty} R_*(t) \geq L \}.
\end{equation}
By Theorem \ref{T3} and (\ref{lemma}) (applied for fixed $L$
to the events $\{ \sigma^t \in {\cal C}_L \}$ for
a sequence of times), we have $P(A_L) \geq 1/2$ for every $L$.
Then, letting $L\to\infty$, we have 
\begin{equation} \label{limsup}
P(\limsup_{t \to \infty} R_*(t) \, = \, \infty) \, \geq \, 1/2.
\end{equation}
It is easy to see that the event $A_{\infty}$ in (\ref{limsup}) occurs 
if and only if its translation,
i.e., the event that $ \limsup_{t \to \infty} d(\partial C_x(t),x) = \infty $,
occurs.
Thus $A_\infty$ is translation-invariant and by the translation-ergodicity
of $P$, (\ref{limsup}) implies that $P(A_\infty) = 1$, as desired.
\fbox{} \\

As mentioned before, in dimension two,
the state $\sigma^t$ does not have a unique limit as $t \to \infty$
\cite{NNS}; thus we are interested in the limits along subsequences of $t$.
Let us introduce the ($\omega$-dependent) set of all limiting states,
\begin{equation}
{\cal W} = {\cal W} (\omega) = \{ \tilde \sigma \in {\cal S} :
\, \exists \, t_k \uparrow \infty \hbox{ so that }
\sigma_x^{t_k} \to \tilde \sigma_x \,\, \forall 
x \in {\bf Z}^2 \}.
%\hbox{ and some } \{ t_k \}_{k \in {\bf N}}
%\hbox{ with } t_k \uparrow \infty \}.
%{\cal W} = {\cal W} (\omega) = \{ \tilde \sigma \in {\cal S} :
%\sigma_x^{t_k(\omega)} (\omega) \to \tilde \sigma_x 
%\, \forall
%x \in {\bf Z}^d \hbox{ and some } \{ t_k(\omega) \}_{k \in {\bf N}}
%\hbox{ with } t_k(\omega) \uparrow \infty \}.
\end{equation}
The following theorem concerns such subsequence limits.
The first statement of the theorem means that there exists
an $\omega$-dependent sequence 
$t'_k \uparrow \infty$
so that $\sigma_o^{t'_k} = (-1)^k$ and $C_o(t_k) \supset Q_{k}$.

\begin{theorem} \label{T4}
$ {\cal W} \supset {\cal C} $ almost surely.
Moreover, 
$ {\cal W}$ contains a non-constant absorbing state with a
flat domain wall passing next to the origin, almost surely.
\end{theorem}

{\bf Proof. } The proof 
that the constant $+1$ (resp., $-1$) state is in $ {\cal W}$ is
%of ${\cal W} \supset {\cal C}$ is 
essentially the same as the proof of
Corollary \ref{C1}, but with the event $A_L$ replaced by the event
that $C_L^+ (t)$ (resp., $C_L^- (t)$) occurs for an unbounded
set of $t$'s. 

To prove the second part of the theorem, we consider the square $Q_L$
(with L even, for a reason to be seen later) and
use the fact \cite{NNS} that $\sigma_o^t$ flips infinitely
many times almost surely. We restrict attention to times
$t$ greater than $T_L$, the almost surely finite time after which
$\Gamma^t (Q_L^*)$ is not e-absent and hence satisfies various
geometric constraints, including those 
discussed in the proof of Lemma \ref{L5}.

There will almost surely be a sequence $t_k \to \infty$, of times
(either just before or just after the flips of the origin) when
there is a (monotonic) domain wall $\gamma$, whose endpoints
are just outside $Q_L^*$, passing next to the origin and 
containing the origin in its spanning rectangle $R(\gamma)$.
There are two possibilities: (1) $\gamma$ connects two opposite
sides of the boundary of $Q_L^*$, or (2) two adjacent sides.

Let $B_L^1$ (resp., $B_L^2$) denote the event that (1) (resp., (2))
occurs for an unbounded set of $t$'s. Then $P(B_L^1)+P(B_L^2) \geq 1$
and so $\liminf_{L\to\infty} P(B_L^i)>0$ either for $i=1$ (we call this
case 1) or $i=2$ (case 2) or both. The $\gamma$ of (1) may be either
horizontal or vertical, so we express $B_L^1$ as the
(not necessarily disjoint) union of events $B_L^{1,h}$ and $B_L^{1,v}$
according to whether a horizontal or a vertical $\gamma$ recurs.
By symmetry, these two events have equal probability, so that
in case 1, there is a subsequence $L_j \to \infty$ such that
$P(B_{L_j}^{1,h})>\delta>0$ for all $j$. 

Consider a time $t>T_{L_j}$
when such a horizontal $\gamma$ is present.
By the monotonicity and other geometric restrictions on the
domain walls that follow because $\Gamma^t(Q_L^*)$ is not e-absent,
it follows that there is a sequence of spin flips (with a bounded away from
zero probability of occurring in the next unit time interval) that
will deform such a $\Gamma^t(Q_L^*)$ into one where there is a
horizontal flat domain wall $\gamma'$ just under the origin
and further such that $\sigma^t \in {\cal A}_L$
(i.e., inside $Q_L$, $\sigma^t$ agrees with an absorbing
state).
Thus the event ${\tilde A}_{L_j}^h$, that {$\sigma^t \in {\cal A}_L$ with
a horizontal flat domain wall just under the origin for an unbounded set of
$t$'s}, has $P({\tilde A}_{L_j}^h) \geq P(B_{L_j}^{1,h})>\delta>0$
for all $j$. 

Proceeding as in the proof of Corollary \ref{C1} (but using ergodicity
only with respect to translations in the first coordinate), we
conclude that $P({\tilde A}_{L}^h \hbox{ occurs for all }L)\,=\,1$,
which, together with standard compactness arguments, completes the proof
for case 1. 

Case 2 is similar, but with an extra ingredient. Here we express 
$B_L^2$ as the union of $B_L^{2,NE}$, where $\gamma$ connects the
North and East sides of the boundary of $Q_L^*$, and the three
other directional possibilities. For a Northeast $\gamma$, we use
spin flips to deform $\Gamma^t(Q_L^*)$ to a $\Gamma'$
containing a contour that is 
horizontal and flat from the origin to the Eastern side of $Q_L^*$,
and such that inside the smaller square $Q_{L/2}(L/2,0)$, centered
at $(L/2,0)$, $\Gamma'$ agrees with an absorbing
state. By translating $(L/2,0)$ to the origin and using translation
invariance, we have in case 2 for some subsequence $L'_j$ that
$P({\tilde A}_{L'_j/2}^h) \geq P(B_{L'_j}^{2,NE})>\delta'>0$
for all $j$. The remainder of the proof is as in case 1. \fbox{}
\\

\begin{remark} \label{R2}
We note that ${\cal W}$ is almost surely strictly larger than
${\cal A}$.
Otherwise, it could not be the case that almost surely every site
flips infinitely many times.
Indeed, by arguments similar to those used for
Theorem \ref{T4}, there almost surely 
must be recurrent states that have a
domain wall passing by the origin that is flat except for a
single step right by (or any fixed distance away from) the origin.
There are of course also many states that are almost surely transient,
such as ones with non-monotonic domain walls or more generally ones
that, restricted to some $Q_L$, do not satisfy the geometric conditions,
such as those
in the proof of Lemma \ref{L5}, that prevent non-monotonic
domain walls from forming.
\end{remark}

\begin{remark} \label{R3}
The proof of Theorem \ref{T4} makes it clear that
${\cal W}$ must almost surely contain non-constant absorbing states
both with horizontal and with vertical flat domain walls. 
Similarly, the single step domain walls mentioned in the previous remark
will be both horizontal and vertical (and with the steps at all
possible distances from the origin). It is a natural
conjecture, discussed in Section 1, that almost surely 
all recurrent states, besides the two constant states,
have only a single (monotonic) doubly-infinite domain wall.
Two possibilities as to the exact class of recurrent states are
discussed at the end of Section~\ref{synopsis}. 
We have not been able to show that almost surely
some state with a single doubly-infinite domain wall is 
in fact recurrent.
However, by using arguments like those in the proof of Theorem \ref{T4},
one can show that there are almost surely recurrent (absorbing)
states with
a flat domain wall next to the origin and no other domain walls in
a half-space. 
\end{remark}

\medskip

\noindent
{\bf Acknowledgments.\/}  The research reported here was supported
in part by the U.S.~NSF under grants DMS-98-02310 (F.C.) and DMS-98-03267
(C.M.N.) and by the Italian CNR (E.D.S.). 
E.D.S.~thanks the Courant Institute for its kind hospitality.

\end{document}